\documentstyle[12pt,amscd]{amsart}
\newtheorem{Theorem}{Theorem}[section]

\newtheorem{Corollary}[Theorem]{Corollary}
\newtheorem{Proposition}[Theorem]{Proposition}

\newtheorem{exm}[Theorem]{Example}
\newtheorem{quest}[Theorem]{Question}

\def\depth{\operatorname{depth}}
\def\reg{\operatorname{reg}}

\def\sk{\smallskip\par}
\def\mm{{\mathfrak m}}
\def\afr{{\mathfrak a}}
\def\bfr{{\mathfrak b}}

\begin{document}
\title{ A note on Castelnuovo-Mumford regularity  and   Hilbert coefficients}
\thanks{Both authors were partially supported by NAFOSTED
(Vietnam), under the grant number 101.04-2015.02. The second author was also supported by the Project VAST.HTQT.NHAT.1/16-18.}

\maketitle

\begin{center} 
LE XUAN DUNG\\
Department of  Natural Sciences,  Hong Duc University\\
307 Le Lai, Thanh Hoa, Vietnam\\
E-mail: lxdung27@@gmail.com\\ [15pt] 
and\\[15pt] 
LE TUAN HOA \\ 
 Institute of Mathematics Hanoi, VAST\\ 18 Hoang Quoc Viet, 10307 Hanoi, Vietnam\\
 E-mail: lthoa@@math.ac.vn
  \end{center}

\begin{abstract}  New upper and lower bounds on the Castelnuovo-Mumford regularity  are given in terms of the Hilbert coefficients. Examples are provided to show that these bounds are in some sense nearly sharp.
\end{abstract}

\date{}
\section{Introduction} \sk

Hilbert coefficients are basic invariants associated to primary ideals and the \\ Castelnuovo-Mumford regularity is one of the most important invariants measuring the complexity of a graded algebra.  It was shown in  \cite[Theorem 17.3.6]{BS}, \cite[Lemma 4]{ST2}  and  \cite[Theorem 2]{Tr1} that one can bound the Castelnuovo-Mumford regularity of a graded algebra in terms of its Hilbert coefficients. These bounds  are recursively defined. In \cite[Lemma 1.2]{DH2}, an explicit  bound was given. However these bounds are far from being sharp. In this note, under an additional assumption, we provide a new upper bound (see Theorem \ref{B}). 
The main meaning of this work is not just to provide another bound, but also to show that the 
new bound  is nearly sharp in any dimension.  

 The approach in \cite{BS} and \cite{ST2} uses induction  on the dimension and an idea of  Mumford; it  works for  graded modules, while our approach, like in the proofs of \cite[Lemma 3.1]{CM} and \cite[Theorem 9]{HH},  uses an extended version of the Gotzmann's regularity theorem  in \cite{Bl}.   Therefore this new bound only holds  for graded algebras.

Using a result in the  Erratum of \cite{DH2} we also provide a rough lower bound on the Castelnuovo-Mumford regularity  of a graded algebra in terms of its Hilbert coefficients, see Proposition \ref{D}. An example is given to show that this bound is also in some sense nearly sharp. 

We also pose a question on bounding  the Castelnuovo-Mumford regularity of the associated graded ring of a local ring in terms of fewer number of Hilbert coefficients and a question on the dependence of Hilbert coefficients.


\section{Results} \label{Re}

Let $R= \oplus_{n\ge 0}R_n$ be a Noetherian standard graded ring over a local Artinian ring $(R_0,\mm_0)$. We  always assume that $R_0/\mm_0$ is an infinite field.  Let $E$ be a finitely generated graded $R$-module of dimension $d$.  
First let us recall some notation. For $0\le i \le d$, put
 $$ a_i(E) =
\sup \{n|\ H_{R_+}^i(E)_n \ne 0 \} ,$$
where $R_+ = \oplus _{n>  0} R_n$. The {\it Castelnuovo-Mumford regularity} of $E$  is defined by 
$$\reg(E) = \max \{ a_i(E) + i \mid  0\le  i \leq  d \},$$
and the {\it Castelnuovo-Mumford regularity of $E$  at and above level $p$}, $0\le p\le d$,  is defined by 
$$\reg^p(E) = \max \{ a_i(E) + i \mid  p\leq i \leq  d \}.$$
We denote the Hilbert function $\ell_{R_0}(E_t)$ and the Hilbert polynomial of $E$ by $h_E(t)$ and $p_E(t)$, respectively. Writing $p_E(t)$ in the form:
$$p_E (t) = \sum_{i=0}^{d-1} (-1)^i e_i(E){t+d-1-i \choose d-1-i},$$
  the numbers $e_i(E)$ are called   {\it Hilbert coefficients} of $E$. For $p\le d-1$, let
$$\xi_{p}(E) = \max\{e_0(E), |e_1(E)|,...,|e_p(E)|\}.$$

It was shown in  \cite[Theorem 17.3.6]{BS}, \cite[Lemma 4]{ST2}  and  \cite[Theorem 2]{Tr1}  (see also  \cite[Lemma 1.2]{DH2} for an explicit formula) that  $\reg^1(E)$ can be  bounded in terms of $\xi_{d-1}(E)$, provided that $E$ is generated in degrees at most 0. In the case of quotient rings of $R$, the following main result provides a much better bound. Its proof  uses an  approach  developed in \cite[Section 3]{HH}.

\begin{Theorem}\label{A}
Let $R_0$ be an Artinian equicharacteristic local ring and $I\subset R_+$ a homogeneous ideal of $R$ such that $\dim R/I  = d\ge 1$.  
Then, for all $1\le p\le d$, we have
$$\reg^p(R/I) \leq (\xi_{d-p}(R/I) +1)^{2^{d-p}}-2.$$
\end{Theorem}

\begin{pf}  We may assume that $I\neq 0$. Let $e_j := e_j(R/I)$ and $\xi_j  := \xi_j (R/I)$.  By \cite[Corollary 3.5(i)]{Bl}, the Hilbert polynomial can be uniquely written in the form
$$p_{R/I}(t)= {c_1+ t \choose t} + {c_2+ t-1 \choose t-1}+ \cdots + {c_s+ t-s+1 \choose t-s+1},$$
where $c_1\geq c_2 \ge \cdots \ge c_s \geq 0$ are integers.  For $0\le j \le d-1$ set
$$B_j  =  B_j(R/I)  = \sharp\{i;\  c_i \ge (d-1)- j\}.$$
Note that $ e_0 = B_0 \le B_1 \le \cdots \le  B_{d-1} =s$ and $\xi_0 \le \xi_1 \le \cdots \le \xi_{d-1}$. By \cite[Corollary 3.5(ii)]{Bl}, $\reg^p(R/I) \le B_{d-p} -1$. Hence, it suffices to show that
  $$B_j \le (\xi_j+1)^{2^j} -1,$$
for all $0 \le j \le d-1$. Since $B_0= e_0 = \xi_0$,  the inequality holds for $j=0$. For $j\ge 1$, by \cite[Proposition 3.9]{Bl} we have
\begin{equation} \label{EA3a}
B_j = (-1)^je_j + {B_{j-1}+1\choose 2} - {B_{j-2}+1\choose 3}+ \cdots + (-1)^{j-1}{B_0+1\choose j+1}.
\end{equation}
For $j=1$ it yields
$$B_1 = -e_1 + {B_0+1 \choose 2} = -e_1 + {e_0+1 \choose 2}\le \xi_1 + {\xi_1+1 \choose 2} < (\xi_1 +1)^2 -1.$$
Let $j\ge 2$. By the induction assumption we may assume that 
$$B_{j-l} \le (\xi _{j-l}+1)^{2^{j-l}} -1 \le  (\xi _{j}+1)^{2^{j-l}} -1,$$
for all $1\le l\le j$. Since $2^l \ge l+1$, we have 
$${B_{j-l} +1 \choose l+1} \le \frac{ ( B_{j-l} +1)^{l+1}}{(l+1)!}  \le  \frac{ ( B_{j-l} +1)^{2^l}}{(l+1)!} \le  \frac{ ( \xi_j+1)^{2^j}}{(l+1)!}.$$
By (\ref{EA3a}) this implies
$$\begin{array}{ll}
B_j  & \le |e_j | + {B_{j-1}+1\choose 2} + {B_{j-3}+1\choose 4}+ \cdots \\
&\le \xi _j+ (\xi_j +1)^{2^j}(\frac{1}{2!} + \frac{1}{4!} + \cdots)\\
&< \xi_j + \frac{7}{12}(\xi_j +1)^{2^j} < (\xi_j +1)^{2^j} -1.
\end{array}$$
This completes the proof of the theorem.
\end{pf}

 Let $I$ be an $\mm$-primary ideal of a $d$-dimensional Noetherian local ring $(A,\mm)$. Then,  for $n \gg 0$,  we can write
$$H_{A/I}(n) : = \ell(A/I^{n+1}) =  \sum_{i=0}^d (-1)^i e_i(I){n+d-i \choose d-i}.$$
The integers $e_i(I)$ are called {\it  Hilbert coefficients} of $I$.  Let $G(I) = A/I \oplus I/I^2 \oplus \cdots $. Note that $e_i(I)=e_i(G(I))$ for $0\le i \le d-1$. 

As an application of the above theorem, we can give a much better bound for $\reg (G(I))$ than the one in \cite[Theorem 1.8]{DH2}. We always assume that $A/\mm$ is infinite.

\begin{Theorem} \label{B}Let $I$ be an $\mm$-primary ideal of a Noetherian local ring $(A, \mm)$ of dimension $d\ge 1$ such that $A/I$  is equicharacteristic , and let
$$\xi_p(I) := \max\{e_0(I), |e_1(I)|,...,|e_p(I)|\}.$$
Then
$$\reg(G(I)) \le  (\xi_d(I) +1)^{2^d}  - 2.$$
Moreover, if $\depth(A) \ge 1$, then 
$$\reg(G(I)) \le  (\xi_{d-1}(I) +1)^{2^{d-1}}  - 2.$$
\end{Theorem}

\begin{pf} 
If $\depth (A) \ge 1$, then by \cite[Theorem 5.2]{H0}, $\reg(G(I)) = \reg^1(G(I))$. Note that $\xi_{d-1} (G(I) )  = \xi_{d-1}(I)$. Hence the second statement follows from Theorem \ref{A}.

In order to prove the first statement, let $x$ be an indeterminate of $\deg(x) =1$ and let $S= G(I)[x]$ be a standard graded ring of dimension $d+1$ over $A/\mm$.  Then 
$h_S(n) = \sum_{i=0}^n \ell_A(I^i/I^{i+1}) = \ell(A/I^{n+1}).$
Hence $e_i(I) = e_i(S)$ for all $i\le d$ and $\xi_d(S) = \xi_d(I)$.  Since $x$ is regular on $S$, $\reg(G(I)) = \reg^1(S)$. Now we can again apply Theorem \ref{A} to $S$ in order to complete the proof.
\end{pf}

As said in the introduction, bounds on $\reg^1(E)$ for any Noetherian graded $R$-module $E$ in terms of Hilbert coefficients were given in \cite{BS} and \cite{ST2}.  A bound on $\reg (G_M(I))$ is given  in \cite[Theorem 1.8]{DH2}, where $I$ is an $\mm$-primary ideal  of a Noetherian local  ring $(A,\mm)$, $M$ is a finite $A$-module and $G_M(I) = \oplus_{n\ge 0}I^nM/I^{n+1}M$. The orders of these bounds are bigger than the one in the above two theorems. We don't know if the bounds $\reg^p(E) \leq (\xi_{d-p}(E) +1)^{2^{d-p}}-2$ and $\reg(G_M(I)) \le  (\xi_{d}(I,M) +1)^{2^d}  - 2$  hold, where $e_i(I,M)$ denotes the $i$-th Hilbert coefficients of $I$ with respect to $M$ and  $\xi_{d}(I,M)  = \max\{e_0(I,M), |e_1(I,M)|, ..., |e_d(I,M)|\}$.

\begin{Corollary} \label{C}With the assumptions in Theorem \ref{B},  we have 
$$(-1)^{i-1}e_i(I) < \frac{7}{12}(\xi_{i-1} (I)+ 1)^{2^i} - e_0(I),$$
 for all $1\le i\le d$,
\end{Corollary}

\begin{pf} Keep the notation in the proof of the previous two theorems. Then by (\ref{EA3a}), we have
$$ (-1)^{i-1}e_i = {B_{i-1}(S)+1\choose 2} - {B_{i-2}(S)+1\choose 3}+ \cdots + (-1)^{i-1}{B_0(S)+1\choose i+1} - B_i(S).$$
In the proof of Theorem \ref{A} we have shown that 
$$B_j(S) \le (\xi_j(S) + 1)^{2^j} - 1 \le (\xi_{i-1}(S)+ 1)^{2^j} - 1 = (\xi_{i-1}+ 1)^{2^j} - 1,$$ 
for all $j\le i-1$.  Since $B_i(S) \ge B_0(S) = e_0$, the last part of computation in the proof of Theorem \ref{A} shows that
$$(-1)^{i-1}e_i < \frac{7}{12}(\xi_{i-1} + 1)^{2^i} - e_0.$$
\end{pf}

It is well-known that $e_1(I) \le {e_0(I) \choose 2}$. Without any assumption on the local ring $A$, one cannot bound $|e_i(I)|$ in terms of $\xi_{i-1}(I)$, see \cite[Example 2.7] {DH2}. On the other hand, when $A$ is a Cohen-Macaulay, both $e_1(I)$ and $e_2(I)$ are non-negative and   $e_2(I) \le {e_1(I)\choose 2}$ (see \cite{KM}). Moreover, in this case, it was first proved in \cite{ST2} that $\reg(G(I))$ and all $|e_i(I)|$, where $1\le i\le d$, are bounded in terms of $e_0(I)$ (see \cite{RTV} for an improvement). Generalizing this fact, it was shown in \cite[Theorem 2.4]{DH2} that 
$\reg(G(I))$ and all $|e_i(I)|$ are bounded in terms of  $\xi_{d-t}(I)$, where $t = \depth A$ and $d-t+1\le i\le d$. However, the bounds in \cite{DH2} are too large. Therefore, in view of Theorem \ref{B} and Corollary \ref{C}, we would like to ask

\begin{quest} {\rm Let $I$ be an $\mm$-primary ideal of a Noetherian local ring $(A, \mm)$ of dimension $d\ge 1$ and depth $t$ such that $A/I$ is  equicharacteristic.
  Do the following inequalities  hold

(i) $ \reg(G(I)) < (\xi_{d-t}(I) +1)^{2^d}$, and 

(ii) $|e_i (I)|< (\xi_{d-t}(I) + 1)^{2^i}$ for all $d-t+1 \le i \le d$?}
\end{quest}

Note, by Theorem \ref{B}, that (i) holds for $t\le 1$. The above question is of interest even in the case of a Cohen-Macaulay ring $A$ and $I = \mm$. The following example shows that the bounds in Theorem \ref{A} and Theorem \ref{B} are almost optimal.

\begin{exm} {\rm Let $f_1,...,f_c$ be a regular sequence of homogeneous  polynomials of degree $\delta \ge \max\{c,\ 36\}$ in $R= K[x_1,...,x_n]$ such that $2d +1< c$, where $d:=n-c\ge 3$. Let $\afr= (f_1,...,f_c)\subset R$.  Then $e_0 := e_0(R/\afr) = \delta^c$ and $\reg(R/\afr) = c\delta - c$. By \cite[Proposition A in  Corrigendum]{DH2} we have 
$$|e_i(R/\afr)| \le e_0(c\delta )^i < e_0\delta^{2d} <  e_0^{1+\varepsilon }/\delta \ \ \text{for \ all } \ 1 \le i \le d-1,$$
where $0< \varepsilon \le 1$ and $\varepsilon \rightarrow 0$ if $c/d \rightarrow \infty $.
Hence $\xi_{d-1} (R/\afr) < e_0^{1+\varepsilon }/\delta $. 

Let $\bfr= \text{lex}(I)$ be the lex-segment ideal of $\afr$, that is the ideal of $R$ generated by all first $h_\afr(m)$ monomials in $R$ of degree $m$ with respect to the lexicographic order, where $m$ runs through all positive integers. This ideal has the same Hilbert function as $\afr$. Hence $p_{R/\bfr}(t) = p_{R/\afr}(t)$, which implies $\xi_{d-1} (R/\bfr)  = \xi_{d-1} (R/\afr) < e_0^{1+\varepsilon }/\delta \le (e_0/4)^{1+\varepsilon }$.  Since $R/\afr$ is a Cohen-Macaulay ring, by \cite[Proposition 12]{HH}, we have
$$\reg^1(R/\bfr)= \reg(R/\bfr)  \ge 9\frac{e_0^{2^{d-1}}}{9^{2^{d-2}}} - 1= 9(e_0/3) ^{2^{d-1}} -1>9 (\xi_{d-1} (R/\bfr))^{(1+ \epsilon )2^{d-2}},$$
where $0< \epsilon = 2/(1+\varepsilon )  - 1<1$ and $\epsilon \rightarrow 1$ when $c/d \rightarrow \infty $. This shows that the bound in Theorem \ref{A} (in the case $p=1$) is almost optimal in the sense that there is no constant $\alpha <1$ such that $\reg^1(R/I) \le \alpha (\xi_{d-1} (R/I))^{ 2^{d-1}}$ for all $I$ and $R$. Note that an upper bound on $\reg(R/\bfr)$ (in terms of $\delta$) in this example was first given in \cite[Lemma 3.1]{CM}. 

For an example in the local case, let $S= K[[x_1,...,x_n]]$ and $A= S/\bfr $, $I = \mm$. Then $G(\mm) \cong R/\bfr$. It was observed in the proof of Theorem \ref{B}, that $e_d(\mm) = e_d((R/\bfr)[x]) = e_d((R/\afr)[x])$, where $x$ is an indeterminate. Again by \cite[Proposition A in  Corrigendum]{DH2} we have
$e_d(\mm) \le e_0(\reg(R/I) +1 )^d  <  e_0^{1+\varepsilon }/\delta $. So, in this case, we still have $\xi_{d} (\mm) < e_0^{1+\varepsilon }/\delta $, and
$$\reg (G(\mm) )= \reg(R/\bfr) > 9 (\xi_d (\mm))^{(1+ \epsilon )2^{d-2}}.$$}
\end{exm}

We can also give  very rough lower bounds for $\reg(R)$ in terms of Hilbert coefficients.
\begin{Proposition}\label{D1}
Let $R_0$ be an Artinian  local ring and $I\subset R_+$ a homogeneous ideal of a standard graded algebra $R= R_0[x_1,...,x_n]$  in $n$ indeterminates such that $\dim R/I  = d\ge 1$. Let  $c= n-d$.  
Then we have
$$e_0(R/I) \le \ell(R_0){\reg(R/I)  +c \choose c}.$$
The equality holds if and only if   $R/I$ is a Cohen-Macaulay ring and its  Hilbert-Poincar\'e series $HP_{R/I}(z) := \sum_{n\ge 0} h_{R/I}(n)z^n$ is equal to
$$\frac{\sum_{i=0}^a \ell(R_0){c+i-1\choose i} z^i}{(1-z)^d},$$
for some $a\ge 0$.
\end{Proposition}

\begin{pf} Put  $a:= \reg(R/I)$.  Without loss of generality, we may assume that $y_1:= x_{c+1},..., y_d:= x_n$ form a filter-regular sequence of $R/I$, that is,  all modules \\
$0:_{R/(I, y_1,...,y_j)R}y_{j+1}$, $0\le j<d$,  are of finite length. Since $\reg(R/(I, y_1,...,y_d)R) \le \reg(R/I) = a$ (see, e.g., \cite[Proposition 18.3.11]{BS}) and $R/(I, y_1,...,y_d)R$ is an epimorphic image of $R_0[x_1,...,x_c]$, we have 
$$ e_0(R/I) \le B:=  \ell_{R_0}( R/(I, y_1,...,y_d)R) \leq \sum_{i=0}^a \ell_{R_0}(R_0[x_1,...,x_c]_i) = \ell(R_0){a +c \choose c}.$$
If $e_0(R/I) =  \ell(R_0){a +c \choose c}$, then $e_0(R/I) = B$ and $\ell_{R_0}( [R/(I, y_1,...,y_d)R]_i)  = {c+i-1\choose i}$ for all $ 0\le i\le a$. This implies that $R/I$ is a Cohen-Macaulay ring and (using, e.g., \cite[Remark 4.1.11]{BH})
$$HP_{R/I}(z) = \frac{\sum_{i=0}^a \ell_{R_0}( [R/(I, y_1,...,y_d)R]_i) z^i}{(1-z)^d} = \frac{\sum_{i=0}^a  \ell(R_0){c+i-1\choose i} z^i}{(1-z)^d}.$$
Conversely, if $R/I$ is a Cohen-Macaulay ring, then  $\reg(R/(I, y_1,...,y_d)R)  = \reg(R/I) = a$, and if also $HP_{R/I}(z) = \frac{\sum_{i=0}^a \ell(R_0){c+i-1\choose i} z^i}{(1-z)^d},$ then from this Hilbert-Poincar\'e series we can compute the Hilbert coefficient (see, e.g., \cite[Proposition 4.1.9]{BH}):
$$e_0(R/I) =  \ell(R_0)\sum_{i=0}^a {c+i-1\choose i} =  {a+c\choose c}.$$
\end{pf}

The lower bound in Proposition \ref{D1} is attained by the ideal $I = (x_1,...,x_c)^{a+1}, \ a\ge 0$.  It is also attained by the so-called Stanley-Reisner ideal of a  cyclic polytope $C(n,d)$, i.e.,  the intersection of all possible monomial ideals generated by $c$ variables in $R= K[x_1,...,x_n]$,  where $c= n-d \ge 1$ (see, e.g., \cite[Subsection 5.2]{BH}). In this case,    $e_0(R/I) = { n \choose c}$ and $\reg(R/I) = d$. 

Below is a lower bound using other Hilbert coefficients. Using the ideal $I = (x_1,...,x_c)^{a+1}, \ a\ge 0$, one can easily see that the bound is nearly sharp.
\begin{Proposition}\label{D}
Let $R_0$ be an Artinian  local ring and $I\subset R_+$ a homogeneous ideal of a standard graded algebra $R= R_0[x_1,...,x_n]$ in $n$ indeterminates such that $\dim R/I  = d\ge 1$. Let  $c= n-d$.  
Then we have
$$\reg (R/I) \geq \max\{ \sqrt[c]{c! e_0(R/I) / \ell(R_0)} - (c+1)/2,\   \sqrt[c+i]{|e_i(R/I)|/\ell(R_0)} - 1 : 1\le i\le d-1 \}.$$
\end{Proposition}

\begin{pf} Keep the notation in the proof of the previous proposition. Then it was shown that $B\le  \ell(R_0){a +c \choose c}.$ 
This implies 
$ e_0 (R/I) \le B \leq \ell(R_0) \frac{[a+ (c+1)/2]^c}{c!},$  and so
$$a \ge  \sqrt[c]{c! e_0(R/I) / \ell(R_0)} - (c+1)/2.$$
Since ${a+c \choose c} \le (a+1)^c$, we also have $B \leq \ell(R_0)(a+1)^c$.
Hence, by \cite[Proposition A in  Corrigendum]{DH2}, for all $1\le i< d$ we get
$$|e_i(R/I)|  \le B(a+1)^i \le \ell(R_0)(a+1)^c(a+1)^i =  \ell(R_0) (a+1)^{c+i},$$
which yields the statement of the proposition.
\end{pf}


\begin{thebibliography}{10}

\bibitem{Bl} C. Blancafort,   Hilbert functions of
graded algebras over Artinian rings,  {\em J. Pure Appl. Algebra}  {\bf
125} (1998) 55 - 78. 

\bibitem{BS} M. P. Brodmann and R. Y. Sharp,   {\em Local cohomology:
an algebraic introduction with geometric applications}. Second edition. (Cambridge
Studies in Advanced Mathematics, 136. Cambridge University Press,
Cambridge, 2013). 

\bibitem{BH} W.  Bruns and J. Herzog,  {\em Cohen-Macaulay rings}. (Cambridge Studies in Advanced Mathematics, 39. Cambridge University Press, Cambridge, 1993).
  
  \bibitem{CM}
M. Chardin, and G. Moreno-Socias,  Regularity of lex-segment ideals: some
  closed formulas and applications, {\em Proc. Amer. Math. Soc.} {\bf 131} (2003) 1093 - 1102.
  
 \bibitem{DH2}
 L. X. Dung and L. T.   Hoa,  Dependence of Hilbert coefficients,  {\em manuscripta math.} {\bf 149} (2016)  235 - 249. Erratum:  {\em manuscripta math.} {\bf 154} (2017) 551 - 552.
 
\bibitem{H0} L. T. Hoa,    Reduction numbers of equimultiple ideals,  {\em J. Pure  Appl. Algebra} {\bf 109} (1996) 
111 - 126.
 
 \bibitem {HH} L. T.  Hoa and E. Hyry,  Castelnuovo-Mumford
regularity of initial ideals,  {\em J. Symb. Comp.} {\bf 38} (2004)  1327 - 1341. 

\bibitem{KM} D. Kirby and H. A.  Mehran,  A note on the coefficients of the Hilbert-Samuel polynomial for a Cohen-Macaulay module,  {\em J. London Math. Soc.} (2) {\bf 25}  (1982) 449 - 457.

\bibitem{ST2} V. Srinivas and V.  Trivedi,    On the Hilbert function of a Cohen-Macaulay local ring,  {\em J. Algebraic Geom.} {\bf 6} (1997) 733 - 751.

\bibitem{Tr1} V.  Trivedi,   Hilbert functions, Castelnuovo-Mumford regularity and uniform Artin-Rees numbers,  {\em manuscripta math.}  {\bf 94} (1997) 485 - 499.
 
 \bibitem{RTV} M. E. Rossi, N. V. Trung and G. Valla,  Castelnuovo-Mumford
regularity and extended degree, {\em Trans. Amer. Math. Soc.} {\bf 355} (2003) 1773 - 1786. 


\end{thebibliography}
\end{document}